\documentclass[11pt]{amsart}
\usepackage{amsmath,amssymb,amsthm,mathtools}
\usepackage[margin=1in]{geometry}
\usepackage{enumitem}
\usepackage[colorlinks=true,linkcolor=blue,citecolor=blue]{hyperref}
\usepackage{comment}
\newtheorem{theorem}{Theorem}[section]
\newtheorem{lemma}[theorem]{Lemma}

\newtheorem{corollary}[theorem]{Corollary}
\theoremstyle{definition}
\newtheorem{definition}[theorem]{Definition}
\newtheorem{example}[theorem]{Example}
\theoremstyle{remark}
\newtheorem{remark}[theorem]{Remark}
\newtheorem{assumption}[theorem]{Assumption}

\numberwithin{equation}{theorem}

\newcommand{\Aut}{\operatorname{Aut}}
\newcommand{\Z}{\mathbb{Z}}
\newcommand{\N}{\mathbb{N}}
\newcommand{\id}{\mathrm{Id}}

\title{Some results on the isotropy groups of non-invertible simple derivations}
\author{Sumit Chandra Mishra, Dibyendu Mondal, Pankaj Shukla, and Sneha Singh}
\address{Indian Institute of Technology Indore, Simrol, Khandwa Road, Indore 453 552, India}
\email{sumitcmishra@gmail.com, sumitcmishra@iiti.ac.in}
\address{Indian Institute of Technology Indore, Simrol, Khandwa Road, Indore 453 552, India}
\email{mdibyendu07@gmail.com, mdibyendu@iiti.ac.in}
\address{Indian Institute of Technology Indore, Simrol, Khandwa Road, Indore 453 552, India}
\email{Pankajshuklashvm23@gmail.com, phd2201141002@iiti.ac.in}
\address{Indian Institute of Technology Indore, Simrol, Khandwa Road, Indore 453 552, India}
\email{mscphd2403141015@iiti.ac.in}
\subjclass[2020]{13N15; 12H05; 13P05}
\keywords{Simple derivations; isotropy groups; polynomial rings}
\date{\today}

\begin{document}

\begin{abstract}
Let $k$ be a field of characteristic zero, let $2\le i<n$, and let $d_i$ be a $k$-derivation of
$k[x_1,\dots,x_i]$. We study the $k$-derivation of $k[x_1,\dots, x_n]$ of the form
$d_n=d_i+\sum_{j=i}^{n-1}\bigl(a_jx_{j+1}+g_j(x_j)\bigr)\partial_{x_{j+1}}$ where $a_j\in k$ and
$g_j\in k[x_j]\setminus k$. Under suitable conditions on $d_i$, we show that the isotropy group of
$d_n$ consists of translations and is either trivial or isomorphic to $(k,+)$. 
\end{abstract}

\maketitle

\section{Introduction}
Let $k$ be a field of characteristic zero. A $k$-derivation on a $k$-algebra $R$ is a $k$-linear map such that $d(xy)=d(x)y+xd(y)$ for any $x,y\in R$. An $d$-ideal of $R$ is an ideal $I$ such that $d(I)\subseteq I$. A $k$-derivation on $R$ is called simple when the only $d$-ideals of $R$ are $(0)$ and $R$. Simple derivations have several applications in commutative and non commutative algebra. Simple derivations can be used in constructing  non-commutative Noetherian simple rings $R[x,d]$; see \cite{Ore33}, in constructing simple Lie rings with product defined as $[a,b]=ad(b)-bd(a)$; see \cite{JJ78}, in constructing non-holonomic irreducible modules over Weyl algebra; see \cite{Cou07}, in characterization of regular rings; see \cite{Sei67}, \cite{Har75}. The isotropy group of a $k$-derivation $d$ of $R$, denoted by $Aut(R)_d$, consists of automorphisms of $R$ which commute with $d$. Studying these groups helps in classifying and understanding affine $k$-algebras. Some studies on isotropy group of locally nilpotent derivations and simple derivations can be seen in \cite{Yan22}, \cite{Yan24}, \cite{MP17}, \cite{BV21}, \cite{BL20}.  A $k$-derivation is called {\it non-invertible} if its image set contains no units. Such derivations allow us to explore skew polynomial rings $R[x,d]$ with the defining relation $rx-xr=d(r)$ where $d(r)$ can never be a unit in $R$. Examples of non-invertible derivations can be seen in \cite{Jor81}, \cite{PS26}, \cite{MMS25'}.

Let $R_n=k[x_1,x_2,\ldots,x_n]$ be the polynomial ring in $n$ variables over a field $k$. Several examples of simple derivations of $R_n$, for $n\geq 2$, have been constructed (see \cite{Now08}, \cite{SP25}, \cite{PS26}, \cite{MMS25'}, \cite{MMS25}, \cite{Kou14}, \cite{Kou12}). However, most of these derivations have images containing units of $R_n$. The study of isotropy groups associated with simple derivations has also received considerable attention recently. In \cite{MP17}, L. G. Mendes and I. Pan proved that the isotropy group of a simple derivation of $R_2$ is trivial. In \cite{BL20}, L. N. Bertoncello and D. Levcovitz proved that the isotropy group of a simple Shamsuddin derivation is trivial. They also proposed the following conjecture: a Shamsuddin derivation $D$ is simple if and only if its isotropy group is trivial. Subsequently, D. Yan proved the converse implication of this conjecture; see \cite[Theorem 2.3]{Yan24}. In 2022, D. Yan proposed a conjecture concerning the isotropy groups of simple derivations of $R_n$, which states that $Aut(R_n)_d$ is conjugate to a subgroup of translation for any simple derivation $d$ of $R_n$. Some examples supporting this conjecture can be seen in \cite{MMS25}, \cite{HY23}. Nevertheless, determining the isotropy groups of simple derivations in higher dimensions remains a challenging problem. Motivated by these developments, in this article, we undertake a systematic study of the isotropy groups associated with certain classes of non-invertible simple derivations of $R_n$, for $n\geq 3$.

Let $i,n$ be positive integers such that $2\leq i<n$. Let $d_i$ be a derivation, defined on $R_i$, such that there does not exist any element $ f\in R_i$ satisfying $d_i(f)=cf+h(x_i)$ where $c\in k$ and $h(x_i)\in k[x_i]$ with $deg_{x_i}\;h(x_i)\leq l$ for some $l\in \mathbb{N}\cup \{\infty\}$. Consider the derivation $d_n:=d_i +\sum_{j=i}^{n-1}(a_j x_{j+1}+g_j(x_j))\partial_{x_{j+1}}$ where $a_j\in k$ and $g_j(x_j)\in k[x_j]\setminus k$ for any $i \leq j \leq n-1$. Then we prove that under suitable conditions on $d_i$, if $Aut(R_i)_{d_i}$ is trivial then isotropy group $Aut(R_n)_{d_n}$ is isomorphic to a subgroup of translations (cf. Theorem \ref{thm:main}). We also show that this derivation $d_n$ is simple if the derivation $d_i$ is simple (cf. Corollary \ref{cor:simple}). Thus, we construct a broader class of simple derivations on $R_n$ which support the conjecture given by D. Yan.

In the final section, we establish a result showing that the simplicity of $d_i$ implies the simplicity of $d_n$, provided that the isotropy group of $d_n$ is contained in the group of translations of $x_n$ and certain suitable conditions on the coefficients $a_j$ and polynomials $g_j$ are satisfied.
\section{Isotropy groups}

Let $k$ be a field of characteristic zero, and let $i,n$ be positive integers such that
$2\le i<n$. Let $d_i$ be a $k$-derivation of the polynomial ring $R_i=k[x_1,\dots,x_i]$, and set
$R_n:=R_i[x_{i+1},\dots,x_n]$. We define the $k$-derivation $d_n$ of $R_n$ by
\begin{equation}\label{eq:dn}
d_n=d_i+\bigl(a_ix_{i+1}+g_i(x_i)\bigr)\partial_{x_{i+1}}+\cdots+\bigl(a_{n-1}x_n+g_{n-1}(x_{n-1})\bigr)\partial_{x_n},
\end{equation}
where $a_j\in k$ and $g_j(x_j)\in k[x_j]\setminus k$ for all $i\le j\le n-1$. Then $d_n|_{R_i}=d_i$.
For $i\le j\le n$, we consider $R_j=k[x_1,\dots,x_j]$ and let $d_j$ be the restriction of $d_n$ to $R_j$;
note that $d_n(R_j)\subseteq R_j$, so $d_n(f)=d_j(f)$ for $f\in R_j$. For $i<j\le n$ we have
$d_n(x_j)=a_{j-1}x_j+g_{j-1}(x_{j-1})$. We write $d$ for $d_n$ when no confusion can arise.

\begin{assumption} 
For $h\in R_n$ and a variable $x_m$, $\deg_{x_m}h$ denotes the degree of $h$ in $x_m$. Let $R$ be a ring. Then for any polynomial $f(X)=\sum_{i=0}^t a_i X^i$ in $R[X]$ with $a_t\neq 0$, we can write $f(X)=\sum_{i\in \mathbb{Z}_{\geq 0}}a_i X^i$ where we set $a_i=0$ for $i>t$. We assume that $\deg_{x_i}g_i(x_i)\le l$, and we set $l_i=l$ and $l_j=\infty$ for $i<j\leq n$.
\end{assumption}

\begin{definition}\label{def:C}
Let $l\in\N\cup\{\infty\}$. We say that $d_i$ satisfies the condition $(C_l)$ if the following holds:
if $d_i(r)=ar+g(x_i)$ for some $r\in R_i$, $a\in k$ and $g(x_i)\in k[x_i]$ with
$\deg_{x_i}g(x_i)\le l$, then $r\in k$ and $g\in k$.
\end{definition}

\begin{remark}
Assuming that $d_i$ satisfies the condition $(C_l)$, we observe that $d_i$ is non-invertible, since taking $a=0$ and $g(x_i)\in k$ shows that its image contains no units.
\end{remark}

\begin{lemma}\label{lem:property}
Suppose $d_i$ satisfies condition $(C_l)$. Then, for every $j\in\{i,\dots,n\}$, the following property $(P_j)$
holds: if $d_n(f)=cf+h(x_j)$ for some $f\in R_j$, $c\in k$ and $h(x_j)\in k[x_j]$ with
$\deg_{x_j}h\le l_j$, then $f\in k$ and $h\in k$.
\end{lemma}

\begin{proof}
The property $(P_i)$ is true by the condition $(C_l)$. Assume that the property $(P_j)$ holds for some $i\le j\le n-1$; we prove the property $(P_{j+1})$.
Let $f\in R_{j+1}$ with $d(f)=cf+h(x_{j+1})$, where $c\in k$ and $h\in k[x_{j+1}]$.
If $f=0$, the result is immediate . Otherwise write $f=\sum_{s=0}^tf_sx_{j+1}^s$ with $f_s\in R_j$ and
$f_t\neq0$. Since $d(x_{j+1})=a_jx_{j+1}+g_j(x_j)$, we have
\begin{equation}\label{eq:expand}
d(f)=\sum_{s=0}^t\bigl(d(f_s)+sa_jf_s\bigr)x_{j+1}^s+g_j(x_j)\sum_{s=1}^tsf_sx_{j+1}^{s-1}.
\end{equation}
Suppose $t\ge1$. Write $h=\sum_{s=0}^r h_s x_{j+1}^s$, where $h_s\in k$ for $0\leq s\leq r$. Using \eqref{eq:expand}, comparing the coefficients of $x_{j+1}^t$ in $d(f)=cf+h$, we get
$d(f_t)=(c-ta_j)f_t+h_t$ with $h_t\in k$, so $f_t\in k^*$ by the property $(P_j)$, and comparing the coefficients of
$x_{j+1}^{t-1}$, we get
\[
d(f_{t-1})=\bigl(c-(t-1)a_j\bigr)f_{t-1}+\bigl(h_{t-1}-tf_t\,g_j(x_j)\bigr).
\]
Here $h_{t-1}-tf_tg_j(x_j)\in k[x_j]$ has degree $\deg g_j\le l_j$ (for $j=i$, this is the assumption $deg\;g_i\leq l$), so the property $(P_j)$ implies that
$tf_tg_j(x_j)\in k$. This contradicts $g_j\notin k$. Hence $t=0$.\\
When $t=0$, We have $d(f)=d(f_0)\in R_j$, so from $d(f)=cf+h$, we have $h\in k$, and the property $(P_j)$ gives $f\in k$.
\end{proof}
Next, we recall the matrix associated with $d_i$ in \cite{MMS25}, in a slightly modified form. Write
\[
d_i=h_1\partial_{x_1}+\cdots+h_i\partial_{x_i},\qquad h_j=d_i(x_j)\in R_i\setminus\{0\}.
\]
For $u,v\in\Z^i_{\ge0}$ we write $v\le u$ if $v_q\le u_q$ for all $q$, and $u\cdot t=\sum_qu_qt_q$.
Let $M$ be an $i\times i$ matrix over $\Z$ and $I$ the identity matrix. We say that $M$ satisfies the
condition $(\ast)$ if the system $(M-I)Y\le0$ has no non-trivial solution $Y\in\Z^i_{\ge0}$.

\begin{definition}\label{def:D}
We say that $d_i$ satisfies condition $(D)$ if for every $1\le j\le i$ the polynomial $h_j$ has a
monomial $\alpha_jx^{u_j}=\alpha_jx_1^{u_{j1}}\cdots x_i^{u_{ji}}$, $\alpha_j\in k^*$, such that every
monomial of $h_j$ divides $x^{u_j}$, and the matrix $A_{d_i}:=[u_{jq}]_{1\le j,q\le i}$ satisfies the condition
$(\ast)$.
\end{definition}

\begin{remark}
If every monomial of $h_j$ divides $x^{u_j}$, then $x^{u_j}$ is the unique monomial of $h_j$ of highest total degree. Condition $(D)$ is satisfied, for instance, when each $h_j$ is a monomial plus a constant and the matrix $A_{d_i}:=[u_{jq}]_{1\le j,q\le i}$ satisfies condition $(\ast)$. Example~\ref{ex:necessity} shows that the divisibility requirement cannot be replaced by the requirement that $x^{u_j}$ merely be a monomial of highest total degree.
\end{remark}
\begin{lemma}\label{lem:restriction}
Suppose $d_i$ satisfies the condition $(D)$, and let $\rho\in\Aut(R_n)_{d_n}$. Then $\rho(R_i)=R_i$, and
$\rho|_{R_i}\in\Aut(R_i)_{d_i}$.
\end{lemma}
\begin{proof}
    Let $m\in\{i+1,\dots,n\}$ and suppose that $\rho(x_q)\in R_m$ for all $1\le q\le i$; this holds for
$m=n$. We show that $\rho(x_q)\in R_{m-1}$ for all $1\le q\le i$. Write
\[
\rho(x_q)=\sum_{s=0}^{t_q}f_{qs}x_m^s,\qquad f_{qs}\in R_{m-1},\ f_{qt_q}\neq0.
\]
Let $t=(t_1,\dots,t_i)$. Since $m-1\ge i$, we have $d(R_{m-1})\subseteq R_{m-1}$, and $d(x_m)=a_{m-1}x_m+g_{m-1}(x_{m-1})$.
Hence
\begin{equation}\label{eq:lhs}
\deg_{x_m}d(\rho(x_q))\le t_q .
\end{equation}
Write $h_q=c_{u_q} x_1^{u_{q1}}x_2^{u_{q2}}\cdots x_i^{u_{qi}}+T$ where $T\in R_i$ with total degree $<u_{q1}+u_{q2}+\cdots +u_{qi}$. For the monomial $c_{u_q}x^{u_q}$ of $h_q$, the coefficient of $x_m ^{u_q \cdot t}$ in the polynomial $\rho(c_{u_q}x^{u_q})=c_{u_q}\prod_{r=1}^i\rho(x_r)^{u_{qr}}$ is $c_{u_q}\prod_{r=1}^i f_{rt_r}^{u_{qr}}$. Since $c_{u_q}\neq 0$ and $f_{rt_r}\neq 0$ for all $1\leq r \leq i$, using the fact that $R_{m-1}$ is an integral domain, we conclude that $c_{u_q}\prod_{r=1}^i f_{rt_r}^{u_{qr}}\neq 0$. Hence $x_m$-degree of $\rho(h_q)$ is $u_q\cdot t$. 
Since $\rho\in\Aut(R_n)_{d_n}$, we have $d(\rho(x_q))=\rho(d(x_q))=\rho(h_q)$.  Combining with \eqref{eq:lhs}, we obtain
\begin{equation}\label{eq:ineq}
t_q\ \ge\ u_{q1}t_1+\cdots+u_{qi}t_i\qquad\text{for all }1\le q\le i,
\end{equation}
that is, $(A_{d_i}-I)t\le0$. By condition $(\ast)$, $t=0$, i.e.\ $\rho(x_q)\in R_{m-1}$ for all $1\le q\le i$.

Repeating this for $m=n,n-1,\dots,i+1$, we get $\rho(R_i)\subseteq R_i$. Since
$\rho^{-1}\in\Aut(R_n)_{d_n}$, the same argument gives $\rho^{-1}(R_i)\subseteq R_i$. Thus
$\rho(R_i)=R_i$, and $\rho|_{R_i}$ commutes with $d_n|_{R_i}=d_i$.
\end{proof}
The following example shows that the divisibility requirement in $(D)$ cannot be dropped, even
if $x^{u_j}$ is chosen among the monomials of highest total degree of $h_j$.

\begin{example}\label{ex:necessity}
Let $i=2$, $n=3$, $w=x_1-x_2$, and $d_2=w^2(\partial_{x_1}+\partial_{x_2})$. Then consider $d_3=d_2+x_2\,\partial_{x_3}$.

Here $h_1=h_2=w^2$. Choosing the monomials $x_1^2$ and $x_2^2$, both of highest total degree, gives
$A_{d_2}=\operatorname{diag}(2,2)$, which satisfies the condition $(\ast)$. However, $x_1x_2$ divides neither $x_1^2$
nor $x_2^2$, so the condition $(D)$ fails for this choice (and for every other choice).

The derivation $d_3$ is locally nilpotent: $d_3(w)=0$, $d_3^2(x_1)=d_3^2(x_2)=0$ and
$d_3^3(x_3)=0$. Put $F=w^2x_3-\tfrac12x_2^2$; then $d_3(F)=w^2x_2-x_2w^2=0$. Hence $Fd_3$ is locally
nilpotent and $[Fd_3,d_3]=-d_3(F)\,d_3=0$, so $\rho=\exp(Fd_3)\in\Aut(R_3)_{d_3}$
(see \cite{Fre17}). Since $(Fd_3)^m=F^md_3^m$,
\[
\rho(x_1)=x_1+w^2F=x_1+w^4x_3-\tfrac12w^2x_2^2\notin R_2 .
\]
Thus the conclusion of Lemma~\ref{lem:restriction} fails. In terms of the proof: $t_1=t_2=1$ and
$f_{11}=f_{21}=w^4$, so the terms of $x_3$-degree $2$ in $\rho(x_1^2)$, $\rho(-2x_1x_2)$ and
$\rho(x_2^2)$ cancel; indeed $\rho(h_1)=\rho(w)^2=w^2$, and \eqref{eq:ineq} would read $1\ge2$.
Note that $d_2$ does not satisfy $(C_l)$, since $d_2(w)=0$.
\end{example}

The next lemma is the key step in showing that every element of the isotropy group is triangular.

\begin{lemma}\label{lem:Ps}
Suppose $d_i$ satisfies the condition $(C_l)$. Let $i\le p<s\le n-1$, $c\in k^*$ and $G\in R_p$. Then there is no
$F\in R_s$ such that
\begin{equation}\label{eq:Ps}
d(F)=a_pF-c\,g_s(x_s)+G .
\end{equation}
\end{lemma}

\begin{proof}
We proceed by induction on $s-p\ge1$. Fix $s$, and suppose $F\in R_s$ satisfies \eqref{eq:Ps}.
Write $g_s(x_s)=\sum_{r=0}^{D}\gamma_rx_s^r$ with $D\ge1$ and $\gamma_D\neq0$. If $F=0$, then
$c\,g_s(x_s)=G\in R_p\subseteq R_{s-1}$, which is impossible since $\deg_{x_s}g_s=D\ge1$. So write
$F=\sum_{r=0}^qe_rx_s^r$ with $e_r\in R_{s-1}$ and $e_q\neq0$. As in \eqref{eq:expand}, the coefficient
of $x_s^r$ in $d(F)$ is
\begin{equation}\label{eq:cof}
 d(e_r)+ra_{s-1}e_r+(r+1)e_{r+1}\,g_{s-1}(x_{s-1}).   
\end{equation}
We use the property $(P_{s-1})$ of
Lemma~\ref{lem:property}; note that all the polynomials $h(x_{s-1})$ to which we apply it below are
constants or of the form $\kappa\,g_{s-1}(x_{s-1})+\kappa'$ with $\kappa,\kappa'\in k$, and so have degree
at most $l_{s-1}$.

\emph{Case $q>D$.} Using \eqref{eq:cof} to express the coefficient of $x_s^q$ in $d(F)$, and comparing the coefficients of $x_s^q$ on both sides of \eqref{eq:Ps}, we obtain $d(e_q)=(a_p-qa_{s-1})e_q$, so
$e_q\in k^*$ by the property $(P_{s-1})$. Since $q-1\ge D\ge1$, comparing the coefficients of $x_s^{q-1}$ in \eqref{eq:Ps} gives
\[
d(e_{q-1})=\bigl(a_p-(q-1)a_{s-1}\bigr)e_{q-1}-\bigl(qe_q\,g_{s-1}(x_{s-1})+c\gamma_{q-1}\bigr),
\]
where $\gamma_{q-1}=0$ if $q-1>D$. By the property $(P_{s-1})$, $qe_q\,g_{s-1}(x_{s-1})\in k$, a contradiction.

\emph{Case $q<D$.} The coefficient of $x_s^D$ on the right-hand side of \eqref{eq:Ps} is
$-c\gamma_D\neq0$, while $\deg_{x_s}d(F)\le q<D$, a contradiction.

\emph{Case $q=D$.} Comparing the coefficients of $x_s^D$ in \eqref{eq:Ps} using \eqref{eq:cof}, we have
$d(e_D)=(a_p-Da_{s-1})e_D-c\gamma_D$, so $e_D\in k^*$ by the property $(P_{s-1})$, and then
\begin{equation}\label{eq:dagger}
(a_p-Da_{s-1})\,e_D=c\gamma_D\neq0 .
\end{equation}
If $D\ge2$, comparing the coefficients of $x_s^{D-1}$ in \eqref{eq:Ps} using \eqref{eq:cof} gives
$d(e_{D-1})=(a_p-(D-1)a_{s-1})e_{D-1}-\bigl(De_D\,g_{s-1}(x_{s-1})+c\gamma_{D-1}\bigr)$, and the property $(P_{s-1})$
gives $De_Dg_{s-1}(x_{s-1})\in k$, a contradiction. Hence $D=1$, and \eqref{eq:dagger} shows
$a_p\neq a_{s-1}$.

If $s=p+1$, then $a_{s-1}=a_p$, a contradiction. This proves the lemma when $s-p=1$. We assume, as the induction hypothesis, that the lemma holds whenever $s-p=r$. We now show the result whenever $s-p=r+1.$

Let $s=p+r+1$, comparing the coefficients of $x_s^0$ in \eqref{eq:Ps} using \eqref{eq:cof} we obtain
\begin{equation}\label{eq:f}
    d(e_0)=a_pe_0-e_1\,g_{s-1}(x_{s-1})+(G-c\gamma_0),
\end{equation}
with $e_0\in R_{s-1}$, $e_1\in k^*$ and $G-c\gamma_0\in R_p$. Thus, \eqref{eq:f} has the same form as \eqref{eq:Ps} with $s$ replaced by $s'=s-1$, and since $s'-p=r$, we get a contradiction from the induction hypothesis. 
\end{proof}

\begin{lemma}\label{lem:triangular}
Suppose $d_i$ satisfies the conditions $(C_l)$ and $(D)$, and let $\rho\in\Aut(R_n)_{d_n}$. Then $\rho(R_j)=R_j$ for
every $i\le j\le n$, and for $i+1\le j\le n$,
\[
\rho(x_j)=\lambda_jx_j+\varphi_j,\qquad \lambda_j\in k^*,\ \varphi_j\in R_{j-1}.
\]
\end{lemma}

\begin{proof}
We use induction on $j$. For $j=i$ the claim is Lemma~\ref{lem:restriction}. Let $i+1\le j\le n$ and
assume $\rho(R_{j-1})=R_{j-1}$. Then $\rho(x_j)\notin R_{j-1}$: otherwise $\rho(x_j)=\rho(r)$ for some
$r\in R_{j-1}$, and injectivity gives $x_j=r\in R_{j-1}$. Let $m\in\{j,\dots,n\}$ be the largest index
such that $x_m$ occurs in $\rho(x_j)$, and write
\[
\rho(x_j)=\sum_{s=0}^tf_sx_m^s,\qquad f_s\in R_{m-1},\ t\ge1,\ f_t\neq0 .
\]
Set $G:=g_{j-1}(\rho(x_{j-1}))$. By the induction hypothesis, $G\in R_{j-1}\subseteq R_{m-1}$. From
$d(\rho(x_j))=\rho(d(x_j))=a_{j-1}\rho(x_j)+G$ and \eqref{eq:expand}, we obtain
\begin{equation}\label{eq:rhoxj}
\sum_{s=0}^t\bigl(d(f_s)+sa_{m-1}f_s\bigr)x_m^s+g_{m-1}(x_{m-1})\sum_{s=1}^tsf_sx_m^{s-1}
=\sum_{s=0}^ta_{j-1}f_sx_m^s+G .
\end{equation}
Comparing the coefficients of $x_m^t$ in \eqref{eq:rhoxj} gives $d(f_t)=(a_{j-1}-ta_{m-1})f_t$, so $f_t\in k^*$ by
Lemma~\ref{lem:property}.

Suppose $t\ge2$. Comparing the coefficients of $x_m^{t-1}$ in \eqref{eq:rhoxj} gives
$d(f_{t-1})=(a_{j-1}-(t-1)a_{m-1})f_{t-1}-tf_t\,g_{m-1}(x_{m-1})$, and Lemma~\ref{lem:property} gives
$tf_tg_{m-1}(x_{m-1})\in k$, a contradiction. Hence $t=1$.

Suppose $m>j$. Comparing the coefficients of $x_m^0$ in \eqref{eq:rhoxj} gives
\[
d(f_0)=a_{j-1}f_0-f_1\,g_{m-1}(x_{m-1})+G,
\]
with $f_0\in R_{m-1}$, $f_1\in k^*$ and $G\in R_{j-1}$. This contradicts Lemma~\ref{lem:Ps} with
$p=j-1$ and $s=m-1$. Hence $m=j$, and $\rho(x_j)=f_1x_j+f_0$ with $f_1\in k^*$ and
$f_0\in R_{j-1}$. Together with $\rho(R_{j-1})=R_{j-1}$, this gives $\rho(R_j)=R_j$.
\end{proof}

We now determine the isotropy group of $d_n$. For $c=(c_{i+1},\dots,c_n)\in k^{n-i}$ we put $c_q:=0$
for $1\le q\le i$, and we denote by $\tau_c$ the translation of $R_n$ given by
\[
\tau_c(x_q)=x_q+c_q\qquad(1\le q\le n).
\]
Clearly $\tau_c\circ\tau_{c'}=\tau_{c+c'}$, so $c\mapsto\tau_c$ is an injective group homomorphism
$(k^{n-i},+)\to\Aut(R_n)$.

\begin{lemma}\label{lem:translation}
Suppose $\deg_{x_i}g_i(x_i)\le l$, and suppose $d_i$ satisfies the conditions $(C_l)$ and $(D)$ and
$\Aut(R_i)_{d_i}=\{\id_{R_i}\}$. Then every $\rho\in\Aut(R_n)_{d_n}$ equals $\tau_c$ for some
$c\in k^{n-i}$.
\end{lemma}

\begin{proof}
Let $\rho\in\Aut(R_n)_{d_n}$. By Lemma~\ref{lem:restriction}, $\rho|_{R_i}\in\Aut(R_i)_{d_i}$, so
$\rho(x_q)=x_q$ for $1\le q\le i$. By Lemma~\ref{lem:triangular}, for $i<j\le n$ we have
\[
\rho(x_j)=\lambda_jx_j+\varphi_j,\qquad \lambda_j\in k^*,\ \varphi_j\in R_{j-1}.
\]
We show by induction on $j\in\{i,\dots,n\}$ that $\rho(x_j)=x_j+c_j$ for some $c_j\in k$.

When $j=i$, $\rho(x_i)=x_i$. Let $i\le j\le n-1$, and suppose that
$\rho(x_j)=x_j+c_j$ with $c_j\in k$. Now, we show that $\rho(x_{j+1})=x_{j+1}+c_{j+1}$. Applying $\rho$ to $d(x_{j+1})=a_jx_{j+1}+g_j(x_j)$, we get
\[
\rho\bigl(d(x_{j+1})\bigr)=a_j\lambda_{j+1}x_{j+1}+a_j\varphi_{j+1}+g_j(x_j+c_j),
\]
and
\[
d\bigl(\rho(x_{j+1})\bigr)=\lambda_{j+1}\bigl(a_jx_{j+1}+g_j(x_j)\bigr)+d(\varphi_{j+1}).
\]
Since $\rho$ commutes with $d$, we have $\rho\bigl(d(x_{j+1})\bigr)=d\bigl(\rho(x_{j+1})\bigr)$, and therefore
\begin{equation}\label{eq:phi}
d(\varphi_{j+1})=a_j\varphi_{j+1}+h(x_j),\qquad\text{where } h(x_j):=g_j(x_j+c_j)-\lambda_{j+1}g_j(x_j)\in k[x_j].
\end{equation}
Here $\varphi_{j+1}\in R_j$ and $\deg h\le\deg g_j\le l_j$ (for $j=i$ this is the assumption
$\deg g_i\le l$). By the property $(P_j)$ of Lemma~\ref{lem:property}, $\varphi_{j+1}\in k$ and
$h\in k$. Let $D=\deg g_j\ge1$, and let $\gamma\in k^*$ be the leading coefficient of $g_j$. The
coefficient of $x_j^D$ in $h$ is $\gamma(1-\lambda_{j+1})$, and it vanishes because $h\in k$. Hence
$\lambda_{j+1}=1$, and $\rho(x_{j+1})=x_{j+1}+c_{j+1}$ with $c_{j+1}:=\varphi_{j+1}\in k$.
\end{proof}

\begin{lemma}\label{lem:shift}
Let $g(x)\in k[x]$ with $D=\deg g\ge1$, and let $c\in k$.
\begin{enumerate}[label=\textup{(\arabic*)}]
\item If $D\ge2$ and $g(x+c)-g(x)\in k$, then $c=0$.
\item If $D=1$ and $g(x)=\gamma x+\beta$, then $g(x+c)-g(x)=\gamma c$.
\end{enumerate}
\end{lemma}

\begin{proof}
(1) Let $\gamma_D$ be the leading coefficient of $g$. The coefficient of $x^{D-1}$ in $g(x+c)-g(x)$
is $D\gamma_Dc$. Since $D-1\ge1$, this coefficient vanishes, and since $\operatorname{char}k=0$, we
get $c=0$. (2) is clear.
\end{proof}

\begin{lemma}\label{lem:constraints}
Let $c\in k^{n-i}$. Then $\tau_c\in\Aut(R_n)_{d_n}$ if and only if
\begin{equation*}
g_j(x_j+c_j)-g_j(x_j)=-a_jc_{j+1}\qquad\text{for all } i\le j\le n-1. \tag{$E_j$}
\end{equation*}
\end{lemma}

\begin{proof}
$\tau_c$ commutes with $d$ if and only if $\tau_c(d(x_q))=d(\tau_c(x_q))$ for $1\le q\le n$. For
$q\le i$, both sides equal $h_q$, since $h_q\in R_i$ is fixed by $\tau_c$ and $d(c_q)=0$. For
$q=j+1$ with $i\le j\le n-1$,
\[
\tau_c\bigl(d(x_{j+1})\bigr)-d\bigl(\tau_c(x_{j+1})\bigr)
=\bigl(a_jx_{j+1}+a_jc_{j+1}+g_j(x_j+c_j)\bigr)-\bigl(a_jx_{j+1}+g_j(x_j)\bigr),
\]
which vanishes if and only if $(E_j)$ holds.
\end{proof}
In the further results, whenever $\deg g_m=1$ for $i\le m\le n-1$, we write
$g_m(x_m)=\gamma_mx_m+\beta_m$, where $\gamma_m\in k^*$ and $\beta_m\in k$.
\begin{lemma}\label{lem:rules}
Let $c\in k^{n-i}$ satisfy $(E_j)$ for all $i\le j\le n-1$, and let $i\le j\le n-1$. Then:
\begin{enumerate}[label=\textup{(\arabic*)}]
\item if $\deg g_j\ge2$, then $c_j=0$;
\item if $a_j=0$, then $c_j=0$;
\item if $c_{j+1}=0$, then $c_j=0$;
\item if $a_j\neq0$ and $c_j=0$, then $c_{j+1}=0$;
\item if $a_j\neq0$ and $\deg g_j=1$, then $c_{j+1}=-\dfrac{\gamma_j}{a_j}\,c_j$.
\end{enumerate}
\end{lemma}

\begin{proof}
(1) The right-hand side of condition $(E_j)$ lies in $k$, so Lemma~\ref{lem:shift}(1) gives $c_j=0$.

(2) and (3) In both cases, $(E_j)$ reads $g_j(x_j+c_j)=g_j(x_j)$. If $\deg g_j\ge2$, then $c_j=0$ by
(1). If $\deg g_j=1$, then Lemma~\ref{lem:shift}(2) gives $\gamma_jc_j=0$, so $c_j=0$.

(4) If $c_j=0$, the left-hand side of $(E_j)$ is $0$, so $a_jc_{j+1}=0$.

(5) By Lemma~\ref{lem:shift}(2), $(E_j)$ reads $\gamma_jc_j=-a_jc_{j+1}$.
\end{proof}

\begin{theorem}\label{thm:main}
Let $d_n$ be as in \eqref{eq:dn}, with $a_j\in k$, and let $\deg_{x_i}g_i(x_i)\le l$ for some
$l\in\N\cup\{\infty\}$. Suppose $d_i$ satisfies the following conditions:
\begin{enumerate}[label=\textup{(\roman*)}]
\item condition $(C_l)$,
\item condition $(D)$, and
\item $\Aut(R_i)_{d_i}=\{\id_{R_i}\}$.
\end{enumerate}
Then
\[
\Aut(R_n)_{d_n}=\{\tau_c:\ c\in k^{n-i}\ \text{satisfies } (E_j) \text{ for all } i\le j\le n-1\}.
\]
More precisely:
\begin{enumerate}[label=\textup{(\alph*)}]
\item If $a_j\neq0$ for all $i\le j\le n-1$, then $\Aut(R_n)_{d_n}=\{\id_{R_n}\}$.
\item Suppose $a_j=0$ for some $j$, and let $j^*:=\max\{j:\ i\le j\le n-1,\ a_j=0\}$.
\begin{enumerate}[label=\textup{(b\arabic*)}]
\item If $\deg g_m=1$ for all $j^*<m\le n-1$, then
\[
\Aut(R_n)_{d_n}=\{\tau_{s\varepsilon}:\ s\in k\}\cong(k,+),
\]
where $\varepsilon=(\varepsilon_{i+1},\dots,\varepsilon_n)\in k^{n-i}$ is defined by
\[
\varepsilon_m=0\ \ (i<m\le j^*),\qquad \varepsilon_{j^*+1}=1,\qquad
\varepsilon_{m+1}=-\frac{\gamma_m}{a_m}\,\varepsilon_m\ \ (j^*<m\le n-1).
\]
\item If $\deg g_m\ge2$ for some $j^*<m\le n-1$, then $\Aut(R_n)_{d_n}=\{\id_{R_n}\}$.
\end{enumerate}
\end{enumerate}
\end{theorem}

\begin{proof}
The first assertion follows from Lemmas~\ref{lem:translation} and~\ref{lem:constraints}. Let
$c\in k^{n-i}$ satisfy $(E_j)$ for all $i\le j\le n-1$, and recall that $c_i=0$.

(a) Since $c_i=0$ and $a_j\neq0$ for all $j$, Lemma~\ref{lem:rules}(4) gives inductively
$c_{i+1}=\cdots=c_n=0$. Hence $\tau_c=\id_{R_n}$.

(b) Since $a_{j^*}=0$, Lemma~\ref{lem:rules}(2) gives $c_{j^*}=0$, and repeated application of
Lemma~\ref{lem:rules}(3) gives $c_j=0$ for all $i\le j\le j^*$. By the maximality of $j^*$, we have
$a_m\neq0$ for all $j^*<m\le n-1$.

(b1) Put $s:=c_{j^*+1}$. Lemma~\ref{lem:rules}(5) gives inductively $c_m=\varepsilon_ms$ for
$j^*<m\le n$. Together with $c_j=0=\varepsilon_js$ for $i<j\le j^*$, this gives $c=s\varepsilon$.
Conversely, for every $s\in k$ the vector $s\varepsilon$ satisfies $(E_j)$ for all $j$:
\begin{itemize}
\item for $i\le j<j^*$, we have $c_j=c_{j+1}=0$;
\item for $j=j^*$, we have $c_{j^*}=0$ and $a_{j^*}=0$;
\item for $j^*<j\le n-1$, $(E_j)$ reads $\gamma_j\varepsilon_js=-a_j\varepsilon_{j+1}s$, which holds by
the definition of $\varepsilon$.
\end{itemize}
Hence $\Aut(R_n)_{d_n}=\{\tau_{s\varepsilon}:s\in k\}$, 
and the map
$s\mapsto\tau_{s\varepsilon}$ is an isomorphism $(k,+)\to\Aut(R_n)_{d_n}$.

(b2) Let $j^*<m\le n-1$ with $\deg g_m\ge2$. By Lemma~\ref{lem:rules}(1), $c_m=0$. By
Lemma~\ref{lem:rules}(3), $c_j=0$ for all $i\le j\le m$. Since $a_j\neq0$ for $m\le j\le n-1$,
Lemma~\ref{lem:rules}(4) gives $c_j=0$ for all $m<j\le n$. Hence $c=0$.
\end{proof}

\begin{remark}\label{rem:special} If $a_j=0$ for all $j$, then $j^*=n-1$ and Theorem~\ref{thm:main}(b1) gives
$\Aut(R_n)_{d_n}=\{(x_1,\dots,x_{n-1},x_n+s):\ s\in k\}$, recovering \cite[Theorem~2.4]{MMS25}
(with the condition (ii) there replaced by the condition $(D)$; cf.\ Example~\ref{ex:necessity}).
\end{remark}

\begin{example}
Let $i=2$, $n=4$, $a_2=0$, $a_3\neq0$, and $g_3(x_3)=\gamma x_3+\beta$ with $\gamma\neq0$. Then $j^*=2$,
and Theorem~\ref{thm:main}(b1) gives $\tau_{s\varepsilon}(x_3)=x_3+s$, $\tau_{s\varepsilon}(x_4)=x_4-\frac{\gamma}{a_3}s$. One checks
directly that $\tau_{s\varepsilon}(d(x_3))=g_2(x_2)=d(x_3)$ and
$\tau_{s\varepsilon}(d(x_4))=a_3\bigl(x_4-\tfrac{\gamma}{a_3}s\bigr)+\gamma(x_3+s)+\beta=d(x_4)$. If instead
$\deg g_3\ge2$, the isotropy group is trivial by Theorem~\ref{thm:main}(b2), although $a_2=0$.
\end{example}

Next, we recall an important result of Shamsuddin \cite{Sha77}; see also \cite[Theorem~13.2.1]{Now94}.

\begin{theorem}[\cite{Sha77}]\label{thm:sham}
Let $R$ be a ring containing $\mathbb{Q}$, and let $d$ be a simple derivation of $R$. Extend $d$ to a
derivation $\tilde d$ of $R[t]$ by setting $\tilde d(t)=at+b$, where $a,b\in R$. Then $\tilde d$ is simple
if and only if there is no $r\in R$ such that $d(r)=ar+b$.
\end{theorem}

\begin{corollary}\label{cor:simple}
Suppose $d_i$ satisfies the condition $(C_l)$ and $\deg_{x_i}g_i(x_i)\le l$. If $d_i$ is simple, then $d_n$ is simple.
\end{corollary}

\begin{proof}
Let $i\le j\le n-1$ and suppose $d_j$ is simple. If $d_j(r)=a_jr+g_j(x_j)$ for some $r\in R_j$, then
the property $(P_j)$ of Lemma~\ref{lem:property} gives $g_j\in k$, a contradiction. By
Theorem~\ref{thm:sham}, $d_{j+1}$ is simple. By induction, $d_n$ is simple.
\end{proof}

In \cite{MP17}, Mendes and Pan proved that the isotropy group of a simple derivation of $R_2$ is
trivial. As a consequence, we have the following.

\begin{corollary}\label{cor:i2}
Let $i=2$, let $d_n$ be as in \eqref{eq:dn} with $a_j\in k$ arbitrary, and let
$\deg_{x_2}g_2(x_2)\le l$. Suppose $d_2$ is simple and satisfies the conditions $(C_l)$ and $(D)$. Then $d_n$ is simple,
and $\Aut(R_n)_{d_n}$ is as described in Theorem~\ref{thm:main}.
\end{corollary}

\begin{proof}
By Corollary~\ref{cor:simple}, $d_n$ is simple. By \cite[Theorem~1]{MP17}, $\Aut(R_2)_{d_2}=\{\id_{R_2}\}$.
Now apply Theorem~\ref{thm:main}.
\end{proof}

We now give examples of such extensions $d_n$ obtained from simple derivations of $R_2$.

\begin{corollary}\label{cor:ex1}
Let $m_1,m_2$ be positive integers such that $m_2\nmid m_1$. For $n\ge3$, let
\[
d_n=(1+x_1^{m_2}x_2)\partial_{x_1}+x_1^{m_1}\partial_{x_2}+(a_2x_3+g_2(x_2))\partial_{x_3}+\cdots+(a_{n-1}x_n+g_{n-1}(x_{n-1}))\partial_{x_n},
\]
where $a_j\in k$ and $g_j(x_j)\in k[x_j]\setminus k$ for $2\le j\le n-1$. Then $d_n$ is simple and
$\Aut(R_n)_{d_n}$ is as described in Theorem~\ref{thm:main}. In particular, $\Aut(R_n)_{d_n}$ is trivial if
all $a_j\neq0$.
\end{corollary}

\begin{proof}
The derivation $d_2=(1+x_1^{m_2}x_2)\partial_{x_1}+x_1^{m_1}\partial_{x_2}$ is simple
\cite[Theorems~4.1 and~6.1]{Kou14}, and it satisfies the condition $(C_\infty)$ by \cite[Proposition 2.2]{PS26}.
For $(D)$, take $u_1=(m_2,1)$ and $u_2=(m_1,0)$: the monomials of $h_1=1+x_1^{m_2}x_2$ divide
$x_1^{m_2}x_2$, and $h_2=x_1^{m_1}$ is a monomial. For $Y=(t_1,t_2)^T\in\Z^2_{\ge0}$, the system
$(A_{d_2}-I)Y\le0$ reads
\[
(m_2-1)t_1+t_2\le0,\qquad m_1t_1-t_2\le0 .
\]
Adding, $(m_1+m_2-1)t_1\le0$, so $t_1=0$, and then $t_2=0$. Hence $(D)$ holds, and the result follows from
Corollary~\ref{cor:i2}.
\end{proof}

\begin{corollary}\label{cor:ex2}
Let $\alpha\ge1$, $\beta\ge1$ and $m\ge\beta+1$ be integers. For $n\ge3$, let
\[
d_n=(1-x_2^{\alpha}x_1^{\beta})\partial_{x_1}+x_1^{m}\partial_{x_2}+(a_2x_3+g_2(x_2))\partial_{x_3}+\cdots+(a_{n-1}x_n+g_{n-1}(x_{n-1}))\partial_{x_n},
\]
where $a_j\in k$, $g_j(x_j)\in k[x_j]\setminus k$ for $2\le j\le n-1$, and $\deg g_2=1$. Then $d_n$ is
simple and $\Aut(R_n)_{d_n}$ is as described in Theorem~\ref{thm:main}.
\end{corollary}

\begin{proof}
The derivation $d_2$ is simple \cite[Theorem~2.6]{SP25}, and it satisfies the condition $(C_1)$ by \cite[Propositions 2.8, 2.9]{PS26}.
For condition $(D)$, take $u_1=(\beta,\alpha)$ and $u_2=(m,0)$. The system $(A_{d_2}-I)Y\le0$ reads
$(\beta-1)t_1+\alpha t_2\le0$ and $mt_1-t_2\le0$; the first gives $t_2=0$, and then the second gives
$t_1=0$. Now apply Corollary~\ref{cor:i2} with $l=1$.
\end{proof}

\section{A Partial Converse}

Let $i,n,R_i,R_n,d_i$ and $d_n$ be as in the preceding section. Let $c\in k$. By $(x_1,\dots,x_{n-1},x_n+c)$, we denote the automorphism of $R_n$ which sends $x_j$ to $x_j$ for all $1\leq j \leq n-1$, and sends $x_n$ to $x_n+c$. In this section, we show that if $d_i$ is
simple and the isotropy group of $d_n$ is a subgroup of $\{(x_1,\dots,x_{n-1},x_n+c):\ c\in k\}$, then, under
suitable conditions on the $a_j$ and $g_j$, the derivation $d_n$ is simple.

For $s\in\{i,\dots,n-1\}$ we set $p(s):=1+\max\{\,j:\ s\le j\le n-1,\ a_j=a_s\,\}\in\{s+1,\dots,n\}.$
\begin{theorem}\label{thm:converse}
Let $d_i$ be a simple $k$-derivation of $R_i$, and let $d_n$ be as in \eqref{eq:dn}. Suppose that for every
$s\in\{i,\dots,n-1\}$, we have $\deg g_m=1$ for all $p(s)\le m\le n-1$. If
\[
\Aut(R_n)_{d_n}\subseteq\{(x_1,\dots,x_{n-1},x_n+c)|\ c\in k\},
\]
then $d_n$ is a simple derivation of $R_n$.
\end{theorem}

\begin{proof}
Suppose that $d_n$ is not simple. Since $d_i$ is simple, there is $s\in\{i,\dots,n-1\}$ such that $d_s$ is
simple and $d_{s+1}$ is not. By Theorem~\ref{thm:sham}, there is $r_s\in R_s$ with
$d_s(r_s)=a_sr_s+g_s(x_s)$. Put $u=x_{s+1}-r_s\in R_{s+1}\setminus k$. Then
\[
d(u)=a_sx_{s+1}+g_s(x_s)-a_sr_s-g_s(x_s)=a_su .
\]
Let $p=p(s)$. By hypothesis, $g_m(x_m)=\gamma_mx_m+\beta_m$ with $\gamma_m\in k^*$ and $\beta_m \in k$ for $p\le m\le n-1$, and
by the maximality in the definition of $p$, $a_m\neq a_s$ for $p\le m\le n-1$. Define $c_p=1$ and
\begin{equation}\label{eq:cm}
c_{m+1}=\frac{\gamma_m}{a_s-a_m}\,c_m\qquad(p\le m\le n-1),
\end{equation}
and let $\rho$ be the endomorphism of $R_n$ given by
\[
\rho(x_q)=x_q\quad(1\le q\le p-1),\qquad \rho(x_q)=x_q+c_qu\quad(p\le q\le n).
\]
We first show that \emph{$\rho$ is an automorphism.} If $p>s+1$, then $u\in R_{s+1}\subseteq R_{p-1}$, so $\rho$ is triangular
with all diagonal coefficients equal to $1$. If $p=s+1$, then $\rho(x_{s+1})=2x_{s+1}-r_s$ and
$\rho(x_q)-x_q\in R_{s+1}$ for $q>s+1$, so $\rho$ is triangular with diagonal coefficients in $\{1,2\}$.
In both cases $\rho$ is an automorphism.

We now show that $\rho$ commutes with $d$. For $q\le p-1$ we have $d(x_q)\in R_{p-1}$, which $\rho$ fixes pointwise,
so $\rho(d(x_q))=d(x_q)=d(\rho(x_q))$. For $q=p$, using $a_{p-1}=a_s$ and $\rho(x_{p-1})=x_{p-1}$,
\[
d(\rho(x_p))=a_{p-1}x_p+g_{p-1}(x_{p-1})+a_sc_pu=a_{p-1}(x_p+c_pu)+g_{p-1}(x_{p-1})=\rho(d(x_p)).
\]
For $p\le m\le n-1$,
\[
d(\rho(x_{m+1}))=a_mx_{m+1}+g_m(x_m)+a_sc_{m+1}u,
\]
\[
\rho(d(x_{m+1}))=a_m(x_{m+1}+c_{m+1}u)+g_m(x_m+c_mu)=a_mx_{m+1}+g_m(x_m)+(a_mc_{m+1}+\gamma_mc_m)u,
\]
and these agree by \eqref{eq:cm}. Hence $\rho\in\Aut(R_n)_{d_n}$.

Thus, $\rho(x_p)-x_p=u\notin k$. 
$\rho$ is not of the form
$(x_1,\dots,x_{n-1},x_n+c)$, which is a contradiction.
\end{proof}

\begin{remark}
The two extreme cases of Theorem~\ref{thm:converse} are:
\begin{enumerate}[label=\textup{(\arabic*)}]
\item $a_i=\cdots=a_{n-1}$ then $p(s)=n$ for every $s$, and no condition on the $g_m$ is needed;
\item $a_j$ are pairwise distinct then $p(s)=s+1$ for any $s\in \{i,\ldots n-1\}$ and in this case, the condition simply becomes that $\deg g_m=1$ for all $i+1\le m\le n-1$.
\end{enumerate}
\end{remark}

\section*{Acknowledgments}
The first-named author would like to thank the Department of Science and Technology of India for
the INSPIRE Faculty fellowship grant IFA23-MA 197 and Anusandhan National Research Foundation (ANRF)
for the ARG-MATRICS Grant (ANRF/ARGM/2025/001651/MTR). The second author would like to thank IIT
Indore for the Young Faculty Research Grant (No.\ IITI/YFRSG/2024-25/Phase-VII/05) and Anusandhan
National Research Foundation (ANRF) for the PMECRG Grant (ANRF/ECRG/2025/003435/PMS). The
third-named author is financially supported by the SRF from CSIR India, Sr.\ No.\
09/1022(16075)/2022-EMR-I. During the preparation of this manuscript, the authors used Claude as an auxiliary
tool  for improving the organization, clarity and language of
the manuscript.



\end{document}